\documentclass[12pt,twoside,a4paper,notitlepage]{article}
\usepackage[latin1]{inputenc}
\usepackage{color}
\usepackage{amsmath}
\usepackage{amsfonts}
\usepackage{amssymb}
\usepackage{graphicx}
\usepackage[normalem]{ulem}

\usepackage[a4paper,margin=2cm,vmargin={1cm,2cm},includeheadfoot]{geometry}

\newtheorem{theorem}{Theorem}[section]

\newcommand{\qed}{\hspace{\stretch{3}}$\square$\vspace{1.8ex}}

\newtheorem{lemma}[theorem]{Lemma}
\newtheorem{cor}[theorem]{Corollary}
\newtheorem{prop}[theorem]{Proposition}

\newtheorem{remark}[theorem]{Remark}

\def\E{\mathbb{E}}

\def\T{\text{T}}

\def\P{\mathbb{P}}

\def\l{\langle}

\def\r{\rangle}
\newcommand{\di}{{\rm d}}

\newcommand{\proof}{{\bf Proof. }}
\newcommand{\ie}{i.e. }

\title{\bf  An Euler-Poisson Scheme for L\'evy driven SDEs}
\author{\small{\sc A. Ferreiro-Castilla$^{\dagger,}$\footnote{Supported by a Royal Society Newton International Fellowship.}  \ \ A.E. Kyprianou\footnote{Department of Mathematical Sciences, University of Bath,
Claverton Down, Bath, BA2 7AY, U.K.} 
\ \ R. Scheichl$^\dagger$%l\footnote{Department of Mathematical Sciences, University of Bath, Claverton Down, Bath, BA2 7AY, U.K.} 
}
}

\begin{document}

\maketitle
\begin{abstract}
We describe an Euler scheme to approximate solutions of L\'evy driven Stochastic Differential Equations (SDE) where the grid points are random and given by the arrival times of a Poisson process. This result extends a previous work of the authors in Ferreiro-Castilla et al. \cite{FKSS12}. We provide a complete numerical analysis of the algorithm to approximate the terminal value of the SDE and proof that the approximation converges in mean square error with rate $\mathcal{O}(n^{-1/2})$. The only requirement  of the methodology is to have exact samples from the resolvent of the L\'evy process driving the SDE; classic examples such as stable processes, subclasses of spectrally one sided L\'evy processes and new families such as meromorphic L\'evy processes (cf. Kuznetsov et al. \cite{KKP2011}) are some examples for which the implementation of our algorithm is straightforward.

\bigskip

\noindent {\sc Key words and phrases}: L\'evy processes, meromorphic L\'evy processes, stochastic differential equations, Euler schemes.
\bigskip

\noindent MSC 2000 subject classifications: 60H10, 65C05 

\end{abstract}

\section{Introduction}\label{Intro}

In many applications such as mathematical finance, insurance mathematics, mathematical biology, physics or engineering there is a need to numerically solve SDEs (see e.g.~\cite{CT,GS72,Situ,S91}). Let $Y:=\{Y_{t}\}_{t\in[0,T]}$ be the solution of the stochastic differential equation 
\begin{equation}\label{SDE_main_intro}
Y_{t}=y_{0}+\int_{0}^{t}a(Y_{s-})\di X_{s}\qquad t\in[0,T]\ ,
\end{equation}
where $a$ is smooth enough so that (\ref{SDE_main_intro}) has a strong solution.
Many studies deal with numerical approximation of (\ref{SDE_main_intro}) when $X:=\{X_{t}\}_{t\in[0,T]}$ is a Wiener process. The complete path of $X$ is numerically intractable
and, ultimately, any numerical scheme can only be based on simulating the increments of the driving process. Therefore, typical approximation schemes rely on Taylor type approximations of the integral. For It\^o integrals with respect to Wiener processes, Taylor expansions of arbitrary order are available and therefore approximations of arbitrary convergence rate (cf. Kloeden and Platen \cite{KP92}). 

Several problems arise when replacing $X$ in (\ref{SDE_main_intro}) by a L\'evy process. For instance, increments of $X$ are not available in general and approximations of the driving process are required. Moreover, multiple stochastic integrals with respect to Poisson measures are more difficult to handle and most numerical schemes are based on modifications of a first order Taylor approximation or an Euler scheme, although higher order schemes can be described as in Baran \cite{Baran09}. The basic Euler scheme for (\ref{SDE_main_intro}) is then
\begin{equation}\label{Euler0}
\widehat{Y}_{0}=y_{0}\ ,\qquad \widehat{Y}_{\frac{(i+1)T}{n}}=\widehat{Y}_{\frac{iT}{n}}+a\left(\widehat{Y}_{\frac{iT}{n}}\right)\left
(X_{\frac{(i+1)T}{n}}-X_{\frac{iT}{n}}\right)\qquad\quad 0\leq i\leq n-1
\end{equation}
where $\{iT/n\}_{0\leq i\leq n}$ is a deterministic partition of $[0,T]$ and $n\in\mathbb{N}$. For the exact Euler scheme, where the increments of the L\'evy process $X$ are available, convergence rates are explicit for the weak and the strong error. The weak error refers to the convergence rate of $|\E[f(Y_{T})]-\E[f(\widehat{Y}_{T})]|$ for a function $f$ in a suitable class. Protter and Talay \cite{PT97} require $f\in\mathcal{C}^{4}(\mathbb{R})$ in addition to some condition on the first moments of $X$ to show $|\E[f(Y_{T})]-\E[f(\widehat{Y}_{T})]|=\mathcal{O}(n^{-1})$. The literature on the strong error estimates is less extensive. The strong error refers to the $p$-th moment, for $p\geq1$, of the pathwise convergence, \ie $\E[\sup_{t\in[0,T]}|Y_{t}-\widehat{Y}_{t}|^{p}]$. It can be inferred from Dereich and Heidenreich \cite{DH11} that under finite second moment of $X$ we also have $\E[\sup_{t\in[0,T]}|Y_{t}-\widehat{Y}_{t}|^{2}]=\mathcal{O}(n^{-1})$. 

The above convergence rates are theoretical since the exact distribution of the increments of L\'evy processes are in general not available and an extra error needs to be incorporated in to the above convergence rates due to the approximation. See for instance Jacod et al. \cite{JKMP05} for a weak error estimate where a fairly general assumption is made for the approximation of the increments of $X$, or Dereich and Heidenreich \cite{DH11} for a strong error estimate where the jump component of $X$ is truncated below a certain threshold. Indeed, the most common approach relies on the L\'evy-It\^o decomposition by removing the jumps below a given threshold transforming the original L\'evy process into a jump diffusion process. Therefore one expects that the final convergence rates depend on the structure of the small jumps. Compound Poisson processes are piecewise constant processes with jumps happening at the arrival times of a Poisson point process. Hence, a more promising modification is to move away from the deterministic equally-spaced grid points in (\ref{Euler0}). A jump-adapted discretization scheme consists in interlacing an equally-spaced grid to approximate the continuous component of the driving process, with a random grid given by the jump times of the purely discontinuous part, as described in Rubenthaler \cite{Rubenthaler03}. In its simplest form, the approximation can perform very poorly when the jump component has paths of infinite $p$-variation with $p$ close to $2$ (recall that all L\'evy processes have finite $2$-variation paths) as shown in Dereich and Heidenreich \cite{DH11}. A more sensible approach is to substitute the small jumps by a Gaussian correction as performed in Dereich \cite{D11}, but this method has its limitations as discussed in Asmussen and Rosi\'nski \cite{AR01}. A novel approach described in Kohatsu-Higa et al. \cite{KOT12} is to approximate the small jumps with an extra compound Poisson process matching a given number of moments of the original driving process provided these moments exist. 
Convergence rates for weak errors are derived under further assumptions on the smoothness of the function $f$. Under the assumption that the L\'evy measure is a regular varying function, the authors in \cite{KOT12} combine this later approach with a high order scheme for the continuous part obtaining arbitrary convergence rates for the weak error.

The aim of this paper is to describe an Euler scheme defined entirely on a random grid, built from the arrival times of a Poisson process. In all the methodologies mentioned above the mesh size - the largest grid step in the Euler approximation - is bounded above by a constant. In our scheme this feature can no longer be assumed as the inter-arrival times of a Poisson process are exponentially distributed. 
%This is the main difficulty in deriving the numerical analysis of our method. 
The origin of this scheme is based in the recent developments of the Wiener-Hopf factorization for L\'evy processes by Kuznetsov \cite{Kuz} and Kuznetsov et al. \cite{KKP2011, KKPvS}. The Wiener-Hopf factorization is a distributional decomposition of the path of a L\'evy process in terms of the running supremum and the running infimum. In Ferreiro-Castilla et al. \cite{FKSS12} this factorization is used to sample from the bivariate distribution of $(X_{t},\sup_{s<t}X_{s})$ by constructing a random walk approximation where the time steps are chosen according to a an exponential distribution, \ie the arrival times of a Poisson process. This scheme effectively constructs a skeleton of the path of $X$ and therefore it is natural to investigate how this skeleton would perform to obtain approximations of (\ref{SDE_main_intro}). Although the skeleton constructs a random walk approximation of the path which captures not only the end point but the supremum over each exponential time step, in the present paper we will consider an Euler scheme for the solution $Y_{T}$ of (\ref{SDE_main_intro}) at the end point only. Therefore, the proposed algorithm is a random modification of the Euler scheme where our assumption is that we can sample exactly from the distribution of $X_{\mathbf{e}(n/T)}$, where $\mathbf{e}(n/T)$ is an exponential distribution with mean $T/n$ independent of $X$. In other words, the grid points in our Euler scheme are given by a Poisson point process with rate $n/T$ denoted by $\mathbf{N}(n/T)$, where the mean $T/n$ plays the role of the mesh size. We will call our scheme the Euler-Poisson scheme. Our analysis does not assume any way of obtaining the distribution of $X_{\mathbf{e}(n/T)}$ and there is no reason why the latter should be easier than the distribution of $X_{1}$, for a general L\'evy process. Nevertheless, for a large class of processes called meromorphic L\'evy processes, see Kuznetsov et al. \cite{KKPvS}, the distribution of $X_{\mathbf{e}(n/T)}$ is explicit and samples from it are available and numerically easy to obtain. Moreover, several popular L\'evy processes in finance can be approximated via the meromorphic class  (cf. \cite{coco,FKSS12,FS11}) and hence the proposed scheme could be taken as an alternative to deal with SDEs driven by such financial models keeping the desired stylized behavior for the driving process. We shall elaborate further on this point later in the paper. It is also worth mentioning here that in contrast to the more classical methods mentioned above, our numerical performance does not depend on the jump structure of $X$ as long as we can sample from $X_{\mathbf{e}(n/T)}$.

The main result of the paper derives the convergence rate for the mean square error of the approximation of $Y_{T}$ showing that $\E[|Y_{T}-\widetilde{Y}_{n}|^{2}]=\mathcal{O}(n^{-1/2})$, where $\widetilde{Y}_{n}$ is the approximation obtained via the Euler-Poisson scheme. The construction of the proposed scheme is based on a random grid and therefore there is no direct way to perform a pathwise numerical analysis. Nevertheless, we will show that when using our scheme to compute $\E[f(Y_{T})]$ for a given function $f$, our methodology is closely related to classical discretization schemes for the Partial Integro-Differential Equation (PIDE) associated with $\E[f(Y_{T})]$.  
The paper is organised as follows. The next section will introduce the basic notation and describe the Euler-Poisson scheme. Section \ref{NumAna} will perform the numerical analysis of our methodology and give the rate of convergence in mean square error. Finally, we collect several remarks and observations about our scheme in Section \ref{remarksEP} regarding feasibility, extensions and its relation with PIDEs. Some of the technical results are collected in an appendix.

\section{The Euler-Poisson scheme}

\subsection{Preliminaries}

Let $(\Omega, \mathcal{F},\{\mathcal{F}_{t}\}_{t\geq0},P)$ be a filtered probability space and let  $Y:=\{Y_{t}\}_{t\in[0,T]}$ be a $\mathbb{R}^{d_{Y}}$-valued adapted stochastic process which is the strong solution of the stochastic differential equation 
\begin{equation}\label{SDE_main}
Y_{t}=y_{0}+\int_{0}^{t}a(Y_{s-})\di X_{s}\qquad t\in[0,T]\ ,
\end{equation}
where $a:=\mathbb{R}^{d_{Y}}\to\mathbb{R}^{d_{Y}}\otimes\mathbb{R}^{d_{X}}$ is a coefficient with smoothness to be specified, 
$X:=\{X_{t}\}_{t\in[0,T]}$ is a $d_{X}$-dimensional square-integrable L\'evy process, $y_{0}\in\mathbb{R}^{d_{Y}}$ and $T<\infty$. Recall that a L\'evy process is a stochastic process issued from the origin which enjoys the properties of having stationary and independent increments with paths that are almost surely right-continuous with left limits. It is a well understood fact that, as a consequence of  this definition,  the law of every L\'{e}vy process is characterised through a triplet $(b,\Sigma,\Pi)$, where $b\in\mathbb{R}^{d_{X}}$, $\Sigma\in\mathbb{R}^{d_{X}\times d_{X}}$  and $\Pi$ is a measure concentrated on $\mathbb{R}^{d_{X}}\backslash\{0\}$ such that $\int_{\mathbb{R}^{d_{X}}}(1\wedge |x|^{2})\Pi(\mathrm{d}x)<\infty$. For square-integrable L\'evy processes we have, for all $t \ge 0$,
\[
\mathbb{E}[{\mathrm e}^{{\rm i}\l\theta, X_{t}\r}]={\mathrm e}^{-t\Psi(\theta)} \qquad \text{for all} \ \ \theta\in\mathbb{R}^{d_{X}},
\]
where
\begin{equation*}
\Psi ( \theta )  \ =\  \mathrm{i}\l b,\theta\r +\frac{1}{2}\l\theta,\Sigma
\Sigma^{\text{T}}\theta
\r+\int_{\mathbb{R}^{d_{X}}}(1-\mathrm{e}^{\mathrm{i}\l\theta, x\r}+\mathrm{i}\l\theta, x\r)\Pi ({\rm d}x)
\end{equation*}
is the so-called characteristic exponent of the process and $\langle\cdot,\cdot\rangle$ is the usual inner product. Furthermore, the L\'evy-It\^o decomposition guaranties that we can decompose $X$ as
\begin{equation}\label{decoX}
X_{t}=\Sigma W_{t}+L_{t}+bt\qquad t\geq0\ ,
\end{equation}
where $W:=\{W_{t}\}_{t\in[0,T]}$ is a $d_{X}$-dimensional Wiener process and $L:=\{L_{t}\}_{t\in[0,T]}$ is a $d_{X}$-dimensional $L^{2}(\Omega, \mathcal{F},P)$ martingale representing the compensated jumps of $X$. To ease the notation in the following derivations we may assume, without lost of generality, that there exists a constant $k\in\mathbb{R}^{+}$ such that
\begin{equation*}
\int_{\mathbb{R}^{d_{X}}}|x|^{2}\Pi(\mathrm{d}x)\leq k^{2}\ ,
\qquad
|\Sigma|\leq k\ ,
\qquad
|b|\leq k
\qquad
\text{ and }
\qquad
|y_0|\leq k
\ .
\end{equation*}
We use indistinctly $|\cdot|$ to denote the Euclidean norm for vectors or the Frobenius norm for matrices.

The following theorem can be proved following Situ \cite{Situ} and sets up the usual assumptions on (\ref{SDE_main}) in order to have a unique strong solution.

\begin{theorem}[{Situ \cite[Section 3.1]{Situ}}]\label{general_sol_SDE}
Consider the SDE driven by a square integrable L\'evy process given in (\ref{SDE_main}). Let $a:=\mathbb{R}^{d_{Y}}\to\mathbb{R}^{d_{Y}}\otimes\mathbb{R}^{d_{X}}$ be a measurable funcion such that
\begin{equation*}%\label{assumption2}
|a(x)-a(x')|\leq k'|x-x'|
\qquad
\text{ and }
\qquad
|a(y_{0})|\leq k'
\end{equation*}
for $x, x'\in\mathbb{R}^{d_{Y}}$ and $k'\in\mathbb{R}^{+}$. Then, equation (\ref{SDE_main}) has a unique strong solution adapted to the filtration generated by $X$, $\mathcal{F}^{X}$, such that
\begin{equation*}%\label{sup_sol}
\E[\sup_{t\in[0,T]}|Y_{t}|^{2}]\leq %K_{1}(1+T+T^{2})e^{K_{1}(1+T)T}=
K_{1}\ ,
\end{equation*}
where $K_{1}$ is a positive constant depending on $k'$ and $T$ only.
\end{theorem}

\begin{remark}\rm
Without loss of generality we can consider $k'=k$ in Theorem \ref{general_sol_SDE}. In the following derivations, constants denoted by $\{K_{i}\}_{i\geq0}$ and $\{\kappa_{i}\}_{i\geq0}$ are constants depending on $k$ and $T$ only which may be renamed without further notice in consecutive equations. Let us introduce here the notation $a\lesssim b$ for two positive quantities $a$ and $b$ such that $a/b$ is uniformly bounded. 
\end{remark}

\subsection{The discretization scheme}\label{EP_scheme_s}

As mentioned in the introduction, this paper is concerned with a modification of the standard Euler scheme which consists in substituting the equally-spaced time steps by exponentially distributed time steps, hence one can think about the grid points in our scheme as the arrival times of a Poisson process. For $n\geq1$, let $\{\mathbf{e}_{i}(n/T)\}_{i\geq1}$ be an i.i.d. sequence of random variables in $(\Omega,\mathcal{F},P)$ where $\mathbf{e}(q)$ denotes an exponential random variable such that $\E[\mathbf{e}(q)]=q^{-1}$ and denote by $\mathcal{G}$ the $\sigma$-algebra generated by $\{\mathbf{e}_{i}(n/T)\}_{i\geq1}$, which we assume independent from $X$; we set $\mathbf{e}_{0}=0$ for convenience. The Euler-Poisson scheme is then given by the discrete Markov chain $\widetilde{Y}:=\{\widetilde{Y}_{t_{i}}\}_{i\geq0}$ defined recursively by
\begin{equation}\label{PEs}
\widetilde{Y}_{t_{i}}:=\widetilde{Y}_{t_{i-1}}+a(\widetilde{Y}_{t_{i-1}})\Delta X_{\mathbf{e}_{i}(n/T)}
\qquad \text{ for }\quad i\geq1 \quad\text{ and }\quad\widetilde{Y}_{0}:=y_{0}\ ,
\end{equation}
where $\Delta X_{\mathbf{e}_{i}}:=X_{\mathbf{e}_{i}(n/T)}-X_{\mathbf{e}_{i-1}(n/T)}\stackrel{d}{=}X_{\mathbf{e}(n/T)}$ and  
\begin{equation*}%\label{random_grid_def}
t_{i}:=\sum_{j=0}^{i}\mathbf{e}_{j}(n/T)\ .
\end{equation*}
Note that $t_{i}\stackrel{d}{=}\mathbf{g}(i,n/T)$, where $\mathbf{g}(\alpha,\beta)$ denotes a Gamma distribution with shape parameter $\alpha$ and rate parameter $\beta$.  We will also denote by $\mathbf{N}(n/T):=\{\mathbf{N}_{t}(n/T)\}_{t\geq0}$ the Poisson process with arrival times $\{t_{i}\}_{i\geq0}$. In the above description the mean $T/n$ is the analog of the mesh size for deterministic spaced Euler schemes. 
%From now on, we might avoid mentioning the dependence of the mesh size to ease the notation if there is no confusion. 
Under the above construction we claim that $\widetilde{Y}_{t_{n}}$ is an approximation of $Y_{T}$ and the task of this paper is to derive the asymptotic behaviour of 
\begin{equation}\label{EP_error1}
\lim_{n\to\infty}\E[|Y_{T}-\widetilde{Y}_{t_{n}}|^{2}]\ .
\end{equation}

Before we proceed with the numerical analysis, let us introduce a new process which stochastically interpolates the Euler-Poisson scheme. Denote by $\iota(t)$ the largest grid point in our scheme before $t$, \ie $\iota(t):=\sup [0,t]\cap\{t_{i}\}_{i\geq0}$, and define
\begin{equation}\label{SDE_int}
\widehat{Y}_{t}:=y_{0}+\int_{0}^{t}a(\widehat{Y}_{\iota(s-)})\di X_{s}
=\widehat{Y}_{\iota(t)}+a(\widehat{Y}_{\iota(t)})(X_{t}-X_{\iota(t)})
\end{equation}
for $ t\in[0,t_{n}\vee T]$.
Notice that for $t\in[t_{i},t_{i+1})$ we have $\widetilde{Y}_{t_{i}}=\widehat{Y}_{t_{i}}=\widehat{Y}_{\iota(t)}$ and hence $\widehat{Y}:=\{\widehat{Y}_{t}\}_{t\in[0,t_{n}\vee T]}$ interpolates, in a random way, the chain $\widetilde{Y}$. Yet another important random variable which is going to play a crucial role in the following derivations is the largest gap of the random grid $\{t_{i}\}_{i\geq0}$ restricted to $[0,T]$. Let us denote this $\mathcal{G}$-measurable random variable by
\begin{equation}\label{def_tau}
\tau:=\sup_{s\in[0,T]}(s-\iota(s))\ .
\end{equation}

\subsection{Main result and feasibility of the Euler-Poisson scheme}\label{main_result_discussion}

Under the above notations we can now formally state the main result of the paper. The result will be proven in Section \ref{NumAna}.

\begin{theorem}\label{TH_main_g}
Under the assumptions of Theorem \ref{general_sol_SDE}, we have
\begin{equation*}
\E[|Y_{T}-\widetilde{Y}_{t_{n}}|^{2}]\leq K_{2} n^{-1/2}\ ,
\end{equation*}
where $K_{2}$ is a positive constant depending on $k$ and $T$ only.
\end{theorem}

It is clear from the preceding section that the Euler-Poisson method is of practical interest only if samples from the distribution of $X_{\mathbf{e}(q)}$ are available. In general, there is no reason why the later distribution is easier to handle than the distribution of $X_{1}$ itself. Nevertheless, recent developments in Wiener-Hopf theory for $1$-dimensional L\'evy processes (see for example Kuznetsov \cite{Kuz} or Kuznetsov et al. \cite{KKP2011, KKPvS}) have provided a rich enough variety of examples for which the necessary distributional sampling can be performed and thus the Euler-Poisson scheme may lead to simpler numerical techniques for (\ref{SDE_main}). This family of processes are named meromorphic L\'evy processes, see Kuznetsov et al. \cite{KKP2011, KKPvS}. The Wiener-Hopf factorization gives a lot more information that what is strictly needed in order to implement the Euler-Poisson scheme, as it involves the running supremum and the running infimum of the process $X$. Let $\overline{X}_t = \sup_{s\leq t} X_s$, $\underline{X}_t: = \inf_{s\leq t}X_s$ and $q>0$. Then, the Wiener-Hopf factorisation states that the random variables $\overline{X}_{\mathbf{e}(q)}$ and $\overline{X}_{\mathbf{e}(q)}- X_{\mathbf{e}(q)}$ are independent. Thanks to the so-called principle of duality, that is to say the equality in law of the pair $\{X_{(t-s)-} - X_t: 0\leq s\leq t\}$ and $\{-X_{s}: 0\leq s\leq t\}$, it follows that $\overline{X}_{\mathbf{e}(q)}- X_{\mathbf{e}(q)}$ is equal in distribution to $-\underline{X}_{\mathbf{e}(q)}$. This leads to the following factorisation of characteristic functions
\begin{equation}\label{WHF}
\mathbb{E}({\rm e}^{{\rm i} \theta X_{\mathbf{e}(q)}}) = \mathbb{E}({\rm e}^{{\rm i}\theta\overline{X}_{\mathbf{e}(q)}})\times\mathbb{E}({\rm e}^{{\rm i}\theta \underline{X}_{\mathbf{e}(q)}}),\qquad \text{for all} \ \ \theta\in\mathbb{R}, 
\end{equation}
known as the Wiener-Hopf factorisation. For the class of meromorphic L\'evy processes, the Wiener-Hopf factors are explicit and hence we can efficiently sample from the distribution of $X_{\mathbf{e}(q)}$ through (\ref{WHF}). Indeed, numerical algorithms involving the computation of $X_{\mathbf{e}(q)}$ for meromorphic L\'evy processes are very easy to implement and robust with respect to the jump structure, see for example Ferreiro-Castilla et al. \cite{FKSS12}. One large subfamily of such processes is the $\beta$-class of L\'evy processes, which also conveniently offers all the desirable properties of better known L\'evy processes that are used in mathematical finance, such as CGMY processes, VG processes or Meixner processes; see for example the discussions in  Ferreiro-Castilla and Schoutens \cite{FS11} and Schoutens and van Damme \cite{SD10}. This discussion may suggest that we would only apply the Euler-Poisson scheme if we are interested in SDEs driven by meromorphic L\'evy processes, but it brings the possibility to study new processes associated to the SDE (\ref{SDE_main}).
For instance, the results in Ferreiro-Castilla et al. \cite{FKSS12} and the ones presented here suggest that we can sample and numerically analyze approximate solutions for SDEs like
\begin{equation*}
Y_{t}=y_{0}+\int_{0}^{t}a(Y_{s-},\overline{X}_{s-})\di X_{s}\qquad
\text{ or }
\qquad
Y_{t}=y_{0}+\int_{0}^{t}a(Y_{s-},\underline{X}_{s-})\di X_{s}\ .
\end{equation*}
To our knowledge, such SDEs have not yet been numerically considered in the literature, but it is not difficult to imagine applications of such processes. For instance, models that appear in stochastic dynamics for population or chemical reactions might be modeled by the above SDEs where the knowledge of $\overline{X}$ can replace the artificial barrier restrictions that are usually imposed on the driving processes due to physical constrains (see e.g.~Situ \cite[Chapter 11]{Situ}).

The Wiener-Hopf factorization is a result for $1$-dimensional L\'evy processes. Nevertheless we can still apply our Euler-Poisson scheme for multidimensional L\'evy processes $X$ provided that they are of the form $\mathbf{A}Z$, where $\mathbf{A}\in\mathbb{R}^{d_{X}}\otimes\mathbb{R}^{d_{Z}}$ and $Z:=\{Z_{t}\}_{t\geq0}$ is a $d_{Z}$-dimensional L\'evy processes with independent components $\{Z_{i}\}_{i=1}^{d_{Z}}$ 
such that $Z_{i}$ belongs to the meromorphic class for $1\leq {i}\leq d_{Z}$. Thanks to the componentwise independence we can perform the Wiener-Hopf factorization in each $\{Z_{i}\}_{i=1}^{d_{Z}}$ and obtain the distribution of $Z_{\mathbf{e}(q)}$ and hence $X_{\mathbf{e}(q)}\stackrel{d}{=}\mathbf{A}Z_{\mathbf{e}(q)}$. Under this construction, the dependence between the components of $X$ is limited to a correlation matrix. On the other hand, multidimensional L\'evy processes are much more difficult to treat numerically than their $1$-dimensional  counterparts and only a limited number of constructions allow efficient numerical simulation. For a comprehensive discussion on the topic we refer to Cont and Tankov \cite[Chapter 5]{CT} and the references therein. Generally speaking, one can always construct a multidimensional L\'evy process with a matrix correlation as described above with a general L\'evy process in their marginals, in which case one confronts the same sampling issues as when dealing with general $1$-dimensional L\'evy processes. A second approach is to perform a univariate time change in a multidimensional Brownian motion, achieving generalized versions of many popular models in finance such as the VG, the NIG or the CIR process - note that such a generalization does not allow one to impose subsets of independent components. For random vectors it is possible to characterize the dependence between components in terms of the marginals by means of a copula, but this technique becomes more involved when applied to stochastic processes (cf. Cont and Tankov \cite[Chapter 5.3]{CT} and Kallsen and Tankov \cite{KT06}) and it is unclear how to numerically analyse such an approach.

\section{Numerical Analysis}\label{NumAna}

The construction of the Euler-Poisson scheme lies on a random grid that is supported on an interval that can be smaller or bigger than $[0,T]$. We will split the mean square error described in (\ref{EP_error1}) between what we denote by the discretization error and the hitting error. To fix ideas, let us write
\begin{equation}\label{EP_error}
|Y_{T}-\widetilde{Y}_{t_{n}}|=|Y_{T}-\widehat{Y}_{t_{n}}|
\leq |Y_{T}-\widehat{Y}_{T}|+|\widehat{Y}_{T}-\widehat{Y}_{t_{n}}|\ ,
\end{equation}
where the first term in the right hand side of the above inequality corresponds to the discretization error and the second term to the hitting error. %We first study the behavior of the stochastic interpolation $\widehat{Y}$ and then derive the convergence rates for the discretization and the hitting error.

\subsection{The discretization error}

Heuristically, the discretization error should behave as the classical Euler scheme for deterministic equally-spaced grid points. In order to see this, we first derive a technical lemma which obtains the analogous result for $\widehat{Y}$ to the one described in Theorem \ref{general_sol_SDE} for $Y$.

\begin{lemma}\label{overline_integrableT}
Under the assumptions of Theorem \ref{general_sol_SDE}, the process $\widehat{Y}$ defined in (\ref{SDE_int}) is adapted to $\mathcal{G}\bigvee\mathcal{F}^{X}$ and such that 
\begin{equation*}
\text{(i) }\ \E[\sup_{t\in[0,T]}|\widehat{Y}_{t}|^{2}]\leq K_{3}\ 
\quad\text{ and }\quad
\text{(ii) }\ \E[\sup_{t\in[0,T]}|\widehat{Y}_{t}|^{2}\;|\;\mathcal{G}]\leq K_{3}\ ,
\end{equation*}
where $K_{3}$ is a positive constant depending on $k$ and $T$ only.
\end{lemma}
\proof
The adaptivity property is clear from the right hand side of (\ref{SDE_int}). The square integrability of the first part (i) follows similarly as in the proof of the well known Theorem \ref{general_sol_SDE}, which we briefly review here for the sake of completeness. Let $\sigma_{N}:=\inf\{t>0\ |\ |\widehat{Y}_{t}|>N\}$, $t\in[0,T]$ and write from the definition of $\widehat{Y}$
\begin{multline}
\frac{1}{3}|\widehat{Y}_{t\wedge\sigma_{N}}|^{2}
\leq
|y_{0}|^{2}
+
\left|\int_{0}^{t\wedge\sigma_{N}}a(\widehat{Y}_{\iota(s)})b\di s\right|^{2}
+
\left|
\int_{0}^{t\wedge\sigma_{N}}a(\widehat{Y}_{\iota(s-)})\di (\Sigma W_{s}+L_{s})\right|^{2}\\
\leq
|y_{0}|^{2}
+
(t\wedge\sigma_{N})k^{2}\int_{0}^{t\wedge\sigma_{N}}|a(\widehat{Y}_{\iota(s)})|^{2}\di s
+
\left|
\int_{0}^{t\wedge\sigma_{N}}a(\widehat{Y}_{\iota(s-)})\di (\Sigma W_{s}+L_{s})\right|^{2}
\ ,\label{stoch_int_lemma}
\end{multline}
where we have used Cauchy-Schwarz inequality for the random Lebesgue integral.
Using the Lipschitz condition of $a$ we further derive the growth condition 
\begin{equation}
|a(x)|^{2}=|a(x)-a(y_{0})+a(y_{0})|^{2}\leq 4k^{2}|x|^{2}+2k^{2}(2k^{2}+1)\leq K_{0}|x|^{2}+K_{0}\ ,\label{ineq_a}\\
\end{equation}
for a constant $K_{0}$ depending on $k$ and $T$ only. Using the later upper bound and the definition of the stopping time $\sigma_{N}$, we conclude that the stochastic integral in (\ref{stoch_int_lemma}) is a square-integrable martingale, to which we apply Doob's inequality and the It\^o isometry to obtain
\begin{eqnarray*}
\frac{1}{3}\E[\sup_{r\leq t\wedge\sigma_{N}}|\widehat{Y}_{r}|^{2}]
&\leq&
k^{2}+tk^{2}
\E\left[\int_{0}^{t\wedge\sigma_{N}}|a(\widehat{Y}_{\iota(s)})|^{2}\di s\right]
+
8k^{2}\E\left[\int_{0}^{t\wedge\sigma_{N}}|a(\widehat{Y}_{\iota(s)})|^{2}\di s\right]
\\
&\leq&
k^{2}+(tk^{2}+8k^{2})
\left(K_{0}\E\left[\int_{0}^{t\wedge\sigma_{N}}|\widehat{Y}_{\iota(s)}|^{2}\di s\right]
+K_{0}t
\right)
\\
&\leq&
\kappa_{1}+\kappa_{1}\int_{0}^{t}\E[\sup_{r\leq s\wedge\sigma_{N}}|\widehat{Y}_{r}|^{2}]\di s\ ,
\end{eqnarray*}
where $\kappa_{1}$ is a constant only depending on $k$ and $T$. Finally, Gronwall's lemma gives
\begin{equation*}
\E[\sup_{r\leq t\wedge\sigma_{N}}|\widehat{Y}_{r}|^{2}]\leq 3\kappa_{1}e^{3\kappa_{1}t}\leq 3\kappa_{1}e^{3\kappa_{1}T}=K_{3}
\end{equation*}
and (i) follows by letting $N\uparrow\infty$.
The second part of the claim follows analogously by noting that $X$ is independent of $\mathcal{G}$; 
therefore, conditioned on $\mathcal{G}$, the stochastic integral
\begin{equation*}
\int_{0}^{t\wedge\sigma_{N}}a(\widehat{Y}_{\iota(s-)})\di (\Sigma W_{s}+L_{s})
\end{equation*}
is a martingale with respect to $\mathcal{F}^{X}$, allowing us to use a conditioned version of Doob's inequality as well as a conditioned version of It\^o isometry. Thus, (ii) follows by the same formal derivations as in (i).
\qed

The following theorem derives the asymptotic behavior for the discretization error which ultimately depends on the random mesh width $\tau$ defined in (\ref{def_tau}). The necessary results to obtain bounds for the moments of $\tau$ are derived in Appendix \ref{Appendix}.

\begin{theorem}\label{unconditional_main}
Under the assumptions of Theorem \ref{general_sol_SDE}, we have
\begin{equation*}
\E[\sup_{t\in[0,T]}|Y_{t}-\widehat{Y}_{t}|^{2}]\leq K_{4}\E[2\tau+\tau^{2}]\lesssim n^{-1}\log(n)\ ,
\end{equation*}
where $K_{4}$ is a positive constant depending on $k$ and $T$ only.
\end{theorem}
\proof Let $t\in[0,T]$ and define
\begin{equation}\label{Z_def}
Z_{t}:=Y_{t}-\widehat{Y}_{t}=
\int_{0}^{t}(a(Y_{s})-a(\widehat{Y}_{\iota(s)}))b\di s+
\int_{0}^{t}(a(Y_{s-})-a(\widehat{Y}_{\iota(s-)}))\di (\Sigma W_{s}+L_{s})\ .
\end{equation}
From Theorem \ref{general_sol_SDE} and Lemma \ref{overline_integrableT}, we deduce that the stochastic integral on the right hand side of (\ref{Z_def}) is a square integrable martingale with respect to the filtration $\mathcal{G}\bigvee\mathcal{F}^{X}$. We apply Cauchy-Schwarz inequality to the random Lebesgue integral and Doob's martingale inequality plus the It\^o isometry to the stochastic integral in (\ref{Z_def}) to end up with
\begin{eqnarray}
\frac{1}{2}\E[\sup_{r<t}|Z_{r}|^{2}]
&\leq& \E\left[\sup_{r<t}\left(\int_{0}^{r}(a(Y_{s})-a(\widehat{Y}_{\iota(s)}))b\di s\right)^{2}\right]\notag\\
&&\qquad\qquad+\ \E\left[\sup_{r<t}\left(\int_{0}^{r}(a(Y_{s-})-a(\widehat{Y}_{\iota(s-)}))\di (\Sigma W_{s}+L_{s})\right)^{2}\right]\notag\\%\label{CS+D}\\
&\leq& k^{2}\E\left[t\int_{0}^{t}|Y_{s}-\widehat{Y}_{\iota(s)}|^{2}\di s\right]
+
8k^{2}\E\left[\int_{0}^{t}|Y_{s}-\widehat{Y}_{\iota(s)}|^{2}\di s\right]\notag\\
&\leq&\kappa_{2}\int_{0}^{t}\E[|Z_{s}|^{2}]+\E[|\widehat{Y}_{s}-\widehat{Y}_{\iota(s)}|^{2}]\di s\notag\\
&\leq&\kappa_{2}\int_{0}^{t}\E[\sup_{r<s}|Z_{r}|^{2}]+\E[|\widehat{Y}_{s}-\widehat{Y}_{\iota(s)}|^{2}]\di s\label{Gronwall_aux}\ ,
\end{eqnarray}
where $\kappa_{2}$ is a positive constant depending on $k$ and $T$ only. The next objective is to use Gronwall's lemma in the later inequality and hence we need to control $|\widehat{Y}_{s}-\widehat{Y}_{\iota(s)}|$. By the fact that $X$ has independent increments and using 
the growth condition of $a(x)$ in (\ref{ineq_a}), we write
\begin{eqnarray}
\E[|\widehat{Y}_{s}-\widehat{Y}_{\iota(s)}|^{2}]
&=&\E[|a(\widehat{Y}_{\iota(s)})|^{2}]\E[|X_{s}-X_{\iota(s)}|^{2}]\notag\\
&\leq&\bigg(K_{0}\E[|\widehat{Y}_{\iota(s)}|^{2}]+K_{0}\bigg)\E[|X_{s}-X_{\iota(s)}|^{2}]\notag\\
&\leq&\bigg(2K_{0}\E[|Z_{\iota(s)}|^{2}]+2K_{0}\E[|Y_{\iota(s)}|^{2}]+K_{0}\bigg)\E[|X_{s}-X_{\iota(s)}|^{2}]\label{gron_aux}\ .
%&\leq&2K_{0}k^{2}(2T+T^{2})\E[|Z_{\iota(s)}|^{2}]+
%K_{0}k^{2}(2K_{1}+1)
%\E[2\tau+\tau^{2}]\label{gron_aux}\ ,
\end{eqnarray}
Now,
%where we used the growth condition of $a(x)$ in (\ref{ineq_a}); the last inequality above comes from 
\begin{equation*}%\label{ineq_X1}
\E[|X_{s}-X_{\iota(s)}|^{2}]\leq k^{2}\E[2\tau+\tau^{2}]\leq k^{2}(2T+T^{2})
\end{equation*}
and so together with (\ref{Gronwall_aux}) and (\ref{gron_aux}), as well as Theorem \ref{general_sol_SDE}, we obtain
\begin{equation*}
\E[\sup_{r<t}|Z_{r}|^{2}]\leq \kappa_{2}\E[2\tau+\tau^{2}]+\kappa_{2}\int_{0}^{t}\E[\sup_{r<s}|Z_{r}|^{2}]\di s
\end{equation*}
where we renamed the constant $\kappa_{2}$ and from where 
\begin{equation*}
\E[\sup_{t\in[0,T]}|Y_{t}-\widehat{Y}_{t}|^{2}]\leq\E[2\tau+\tau^{2}]\kappa_{2}e^{T\kappa_{2}}
=\E[2\tau+\tau^{2}]K_{4}
\end{equation*}
follows by Gronwall's inequality. The claim of the theorem is proved up to the computation of $\E[2\tau+\tau^{2}]$ which is given by Proposition \ref{tau_Poisson} in the Appendix \ref{Appendix}.
\qed
 
\subsection{The hitting error}

The next result derives the asymptotic behaviour for the hitting error. We will see that the error boils down to measure how fast the random time $t_{n}$ converges to $T$, which in turn is controlled by the variance of a Gamma distribution. Before we proceed, let us first derive two technical lemmas in the same spirit as Lemma \ref{overline_integrableT}.

\begin{lemma}\label{overline_integrableG_aux}
Under the assumptions of Theorem \ref{general_sol_SDE}, the process $\widehat{Y}$ defined in (\ref{SDE_int}) is adapted to $\mathcal{G}\bigvee\mathcal{F}^{X}$ and such that 
\begin{equation*}%\label{statementW}
\max_{0\leq i\leq n}\E[|\widehat{Y}_{t_{i}}|^{2}]\leq K_{5}\ ,
\end{equation*}
where $K_{5}$ is a positive constant depending on $k$ and $T$ only.
\end{lemma}
\proof Fix $i>0$ and recall the definition of $\widehat{Y}_{t_{i}}$ in (\ref{SDE_int}) to write 
\begin{multline}\label{ite_0}
\E[|\widehat{Y}_{t_{i}}|^{2}]
=\E[|\widehat{Y}_{t_{i-1}}|^{2}]+\E[|a(\widehat{Y}_{t_{i-1}})|^{2}]\E[|X_{t_{i}}-X_{t_{i-1}}|^{2}]+2\E[\widehat{Y}_{t_{i-1}}^{\T}a(\widehat{Y}_{t_{i-1}})]\E[X_{t_{i}}-X_{t_{i-1}}]\\%\label{ite_0}\\
\leq\E[|\widehat{Y}_{t_{i-1}}|^{2}]\left(1+K_{0}2k^{2}\frac{T}{n}\left(1+\frac{T}{n}\right)+2\sqrt{K_{0}}k\frac{T}{n}\right)+K_{0}2k^{2}\frac{T}{n}\left(1+\frac{T}{n}\right)+2\sqrt{K_{0}}k\frac{T}{n}\ ,
\end{multline}
where we used that $t_{i}-t_{i-1}\stackrel{d}{=}\mathbf{e}(n/T)$ and the orthogonal decomposition of $X$ in (\ref{decoX}); we also used the 
%\begin{equation}%\label{ineq_X}
%|\E[X_{t_{i}}-X_{t_{i-1}}]|\leq k\frac{T}{n}\ ,\qquad
%\E[|X_{t_{i}}-X_{t_{i-1}}|^{2}]\leq2k^{2}\frac{T}{n}\left(1+\frac{T}{n}\right)\ ;
%\end{equation} 
growth condition (\ref{ineq_a}) and the following upper bound, which is derived from the assumptions on $a(x)$:
\begin{equation*}
|x^{\T}a(x)|\leq \sqrt{K_{0}}|x|^{2}+\sqrt{K_{0}}\ .
\end{equation*}
It is then clear from (\ref{ite_0}) that there exists a constant $\kappa_{3}$, depending on $k$ and $T$ only, such that
\begin{equation*}
\E[|\widehat{Y}_{t_{i}}|^{2}]
\leq\E[|\widehat{Y}_{t_{i-1}}|^{2}]\left(1+\frac{\kappa_{3}}{n}\right)+\frac{\kappa_{3}}{n}
\leq|y_{0}|^{2}\left(1+\frac{\kappa_{3}}{n}\right)^{i}+i\exp\left(i\frac{\kappa_{3}}{n}\right)\frac{\kappa_{3}}{n}\ ,%\label{ite_i}
\end{equation*}
which follows from the argument that if $x_{m+1}\leq \alpha x_{m}+\beta$ and $\alpha\geq1$, then $x_{m}\leq \alpha^{m}x_{0}+me^{m(\alpha-1)}\beta$. Finally,
\begin{equation*}%\label{ii_first}
\max_{0\leq i\leq n}\E[|\widehat{Y}_{t_{i}}|^{2}]\leq
|y_{0}|^{2}\left(1+\frac{\kappa_{3}}{n}\right)^{n}+e^{\kappa_{3}}\kappa_{3}
\leq
e^{\kappa_{3}}\left(k^{2}+\kappa_{3}\right)
\end{equation*}
which concludes the proof. \qed

\begin{lemma}\label{overline_integrableG}
Under the assumptions of Theorem \ref{general_sol_SDE}, the process $\widehat{Y}$ defined in (\ref{SDE_int}) is adapted to $\mathcal{G}\bigvee\mathcal{F}^{X}$ and such that 
\begin{equation*}
\text{(i) }\ \E[\max_{0\leq i\leq n}|\widehat{Y}_{t_{i}}|^{2}]\leq K_{6}\ 
\quad\text{ and }\quad
\text{(ii) }\ \E\left[\left(\E\left[\max_{0\leq i\leq n}|\widehat{Y}_{t_{i}}|^{2}|\mathcal{G}\right]\right)^{2}\right]\leq K_{6}\ ,
\end{equation*}
where $K_{6}$ is a positive constant depending on $k$ and $T$ only.
\end{lemma}
\proof Define $\Delta \widehat{Y}_{i}:=\widehat{Y}_{t_{i+1}}-\widehat{Y}_{t_{i}}$ and we use the same principles as in (\ref{ite_0}) and Lemma \ref{overline_integrableG_aux} to derive 
\begin{equation*}
\E[|\Delta\widehat{Y}_{t_{i}}|^{2}]
=\E[|a(\widehat{Y}_{t_{i}})|^{2}]\E[|X_{t_{i+1}}-X_{t_{i}}|^{2}]
\leq 
(K_{0}K_{5}+K_{0})
2k^{2}\frac{T}{n}\left(1+\frac{T}{n}\right)
\end{equation*}
for $0\leq i\leq n-1$ and hence
\begin{equation}\label{ii_first_aux}
\max_{0\leq i\leq n-1}\E[|\Delta\widehat{Y}_{t_{i}}|^{2}]\leq \kappa_{4}\frac{1}{n}\ ,
\end{equation}
for some constant $\kappa_{4}$ depending on $k$ and $T$ only. Consider now the filtration $\mathcal{H}_{i}:=\sigma\langle\widehat{Y}_{t_{j}},\ 0\leq j\leq i\rangle$ and the auxiliary random variables
\begin{equation*}
Z_{i}:=\Delta\widehat{Y}_{t_{i}}-\E[\Delta\widehat{Y}_{t_{i}}|\mathcal{H}_{i}]
\end{equation*}
for $0\leq i\leq n-1$. It is clear that $Z_{i}$ is $\mathcal{H}_{i+1}$-measurable and it is not difficult to check that $\sum_{j=0}^{i}Z_{i}$ is a martingale such that $\E[Z_{i}Z_{j}]=0$ if $i\neq j$. Therefore we can write 
\begin{eqnarray}
\max_{0\leq i\leq n}|\widehat{Y}_{t_{i}}|^{2}
&\leq&2\left(|y_{0}|^{2}+\max_{0\leq i\leq n-1}\left|\sum_{j=0}^{i}\Delta\widehat{Y}_{t_{j}}\right|^{2}\right)\notag\\
&=&2\left(|y_{0}|^{2}+\max_{0\leq i\leq n-1}\left|\sum_{j=0}^{i}Z_{j}+\E[\Delta\widehat{Y}_{t_{j}}|\mathcal{H}_{j}]\right|^{2}\right)\notag\\
&\leq&2\left(|y_{0}|^{2}+2\underbrace{\max_{0\leq i\leq n-1}\left|\sum_{j=0}^{i}Z_{j}\right|^{2}}_{(*)}
+2\underbrace{\max_{0\leq i\leq n-1}\left|\sum_{j=0}^{i}\E[\Delta\widehat{Y}_{t_{j}}|\mathcal{H}_{j}]\right|^{2}}_{(**)}\right)\ .\label{auxZ1}
\end{eqnarray}
We now use Doob's martingale inequality and the orthogonality of $\{Z_{i}\}_{i=0}^{n-1}$ to bound $(*)$. Combining this with Jensen's inequality and (\ref{ii_first_aux}) we find that 
\begin{equation*}
\E[(*)]
\leq \E\left[\sum_{j=0}^{n-1}|Z_{j}|^{2}\right]
\leq 2\E\left[\sum_{j=0}^{n-1}|\Delta\widehat{Y}_{t_{j}}|^{2}+|\E[\Delta\widehat{Y}_{t_{j}}|\mathcal{H}_{j}]|^{2}\right]
\leq
4\sum_{j=0}^{n-1}\E\left[|\Delta\widehat{Y}_{t_{j}}|^{2}\right]\leq 4\kappa_{4}\ .
\end{equation*}
Similarly one obtains
\begin{equation*}
\E[(**)]
\leq\E\left[\left(\sum_{j=0}^{n-1}|\E[\Delta\widehat{Y}_{t_{j}}|\mathcal{H}_{j}]|\right)^{2}\right]
=\E\left[\left(\sum_{j=0}^{n-1}
|a(\widehat{Y}_{t_{j}})|k\frac{T}{n}
\right)^{2}\right]\leq k^{2}T^{2}(K_{0}K_{5}+K_{0})%\label{starZ2}
\end{equation*}
using Lemma \ref{overline_integrableG_aux}. The first part (i) follows 
substituting the upper bounds $\E[(*)]$ and $\E[(**)]$ into (\ref{auxZ1}).

For the second part (ii) we consider $\mathcal{H}_{i}:=\mathcal{G}\bigvee\sigma\langle\widehat{Y}_{t_{j}},\ 0\leq j\leq i\rangle$ and reproduce the above derivations up to (\ref{auxZ1}). By the definition of $\Delta\widehat{Y}_{t_{i}}$, we have
\begin{eqnarray}
\E[(*)|\mathcal{G}]
&\leq& 4\sum_{j=0}^{n-1}\E\left[|\Delta\widehat{Y}_{t_{j}}|^{2}|\mathcal{G}\right]
\leq 8k^{2} \max_{0\leq i\leq n-1}\E\left[|a(\widehat{Y}_{t_{i}})|^{2}|\mathcal{G}\right]\sum_{j=0}^{n}\mathbf{e}_{j}\left(1+\mathbf{e}_{j}\right)\label{condA1}\\
\E[(**)|\mathcal{G}]
&\leq& \E\left[\left.\left(\sum_{j=0}^{n-1}
|a(\widehat{Y}_{t_{j}})|k\mathbf{e}_{j}
\right)^{2}\right|\mathcal{G}\right]
\leq k^{2}\max_{0\leq i\leq n-1}\E\left[|a(\widehat{Y}_{t_{i}})|^{2}|\mathcal{G}\right]n\sum_{j=0}^{n}\mathbf{e}_{j}^{2}\ .
\label{condA2}
\end{eqnarray}
Therefore, to prove (ii) it is enough to recall (\ref{auxZ1}) and then show that
\begin{equation*}%\label{fin_aux}
\E\left[(\E[(*)|\mathcal{G}])^{2}+(\E[(**)|\mathcal{G}])^{2}\right]\leq \kappa_{5}
\end{equation*}
for some constant $\kappa_{5}$ depending on $k$ and $T$ only. Using Cauchy-Schwarz inequality, 
a sufficient condition for this claim to be true is
\begin{eqnarray}
\E\left[\left(\max_{0\leq i\leq n-1}\E\left[|a(\widehat{Y}_{t_{i}})|^{2}|\mathcal{G}\right]\right)^{4}\right]&\leq& \kappa_{5}\label{final_aux_i}\\
\E\left[\left(n\sum_{j=0}^{n}\mathbf{e}_{j}^{2}\right)^{4}+
\left(\sum_{j=0}^{n}\mathbf{e}_{j}\left(1+\mathbf{e}_{j}\right)\right)^{4}\right]&\leq& \kappa_{5}
\label{final_aux_ii}
\end{eqnarray}
by renaming $\kappa_{5}$. Since $\E[\mathbf{e}(n/T)^{i}]=i!T^{i}/n^{i}$ for $i\geq1$, one can check that
\begin{equation}\label{first_exp}
\E\left[\left(n\sum_{j=0}^{n}\mathbf{e}_{j}^{2}\right)^{4}+
\left(\sum_{j=0}^{n}\mathbf{e}_{j}\left(1+\mathbf{e}_{j}\right)\right)^{4}\right]\leq 
8!T^{8}+8\left(4!T^{4}+\frac{8!T^{8}}{n^{4}}\right)
\end{equation}
and (\ref{final_aux_ii}) holds. Adapting the left hand side of (\ref{ite_0}) to incorporate the conditional expectation we write
\begin{equation*}
\E[|\widehat{Y}_{t_{i}}|^{2}|\mathcal{G}]
\leq\E[|\widehat{Y}_{t_{i-1}}|^{2}|\mathcal{G}]
\left(1+K_{0}2k^{2}(\mathbf{e}_{i}+\mathbf{e}_{i}^{2})+2\sqrt{K_{0}}k\mathbf{e}_{i}\right)+K_{0}2k^{2}(\mathbf{e}_{i}+\mathbf{e}_{i}^{2})+2\sqrt{K_{0}}k\mathbf{e}_{i}\ ,
\end{equation*}
%\begin{eqnarray*}
%\E[|\widehat{Y}_{t_{i}}|^{2}|\mathcal{G}]
%&\leq&\E[|\widehat{Y}_{t_{i-1}}|^{2}|\mathcal{G}]
%\left(1+K_{0}2k^{2}(\mathbf{e}_{i}+\mathbf{e}_{i}^{2})+2\sqrt{K_{0}}k\mathbf{e}_{i}\right)+K_{0}2k^{2}(\mathbf{e}_{i}+\mathbf{e}_{i}^{2})+2\sqrt{K_{0}}k\mathbf{e}_{i}\\
%&\leq&\E[|\widehat{Y}_{t_{i-1}}|^{2}|\mathcal{G}]\bigg(1+\kappa_{6}(\mathbf{e}_{i}+\mathbf{e}_{i}^{2})\bigg)+\kappa_{6}(\mathbf{e}_{i}+\mathbf{e}_{i}^{2})\ ,
%\end{eqnarray*}
where $\kappa_{6}$ as a constant only depending on $k$. Using again a recurrence argument we easily see that 
\begin{equation*}
\max_{0\leq i\leq n}\E[|\widehat{Y}_{t_{i}}|^{2}|\mathcal{G}]
\leq|y_{0}|^{2}\prod_{i=1}^{n}\bigg(1+\kappa_{6}(\mathbf{e}_{i}+\mathbf{e}_{i}^{2})\bigg)+\sum_{i=1}^{n}\kappa_{6}(\mathbf{e}_{i}+\mathbf{e}_{i}^{2})\prod_{j=i+1}^{n}\bigg(1+\kappa_{6}(\mathbf{e}_{j}+\mathbf{e}_{j}^{2})\bigg)\ .
\end{equation*}
Finally (\ref{final_aux_i}) follows from the same sort of manipulations used in (\ref{first_exp}) and the above inequality. Since (\ref{final_aux_i}) and (\ref{final_aux_ii}) hold, so do (\ref{condA1}) and (\ref{condA2}), which proves (ii). \qed

\begin{prop}\label{propo_hitting}
Under the assumptions of Theorem \ref{general_sol_SDE}, we have
\begin{equation*}
\E[|\widehat{Y}_{T}-\widehat{Y}_{t_{n}}|^{2}]\leq K_{7} n^{-1/2}\ .
\end{equation*}
where $K_{7}$ is a positive constant depending on $k$ and $T$ only.
\end{prop}
\proof
Let us write
\begin{equation*}%\label{decomp_hitting}
\frac{1}{2}|\widehat{Y}_{T}-\widehat{Y}_{t_{n}}|^{2}
\leq\left|
\int_{t_{n}}^{T}a(\widehat{Y}_{\iota(s)})b\di s\right|^{2}
+
\left|
\int_{t_{n}}^{T}a(\widehat{Y}_{\iota(s-)})\di (\Sigma W_{s}+L_{s})\right|^{2}\ .
\end{equation*}
According to part (i) of Lemmas \ref{overline_integrableT} and \ref{overline_integrableG}, the stochastic integral in the above decomposition is a square integrable martingale with respect to $\mathcal{G}\bigvee\mathcal{F}^{X}$; hence, we can use again Cauchy-Schwarz inequality to the random Lebesgue integral and the It\^o isometry for the stochastic integral to obtain 
\begin{eqnarray}
\frac{1}{2}\E[|\widehat{Y}_{T}-\widehat{Y}_{t_{n}}|^{2}]
&\leq&
\E\left[\left(k^{2}|T-t_{n}|+2k^{2}\right)\int_{t_{n}}^{T}|a(\widehat{Y}_{\iota(s)})|^{2}\di s\right]\notag\\
&&\hspace{-3cm}=\ 
k^{2}\E\left[\left(|T-t_{n}|+2\right)\int_{t_{n}}^{T}\E[|a(\widehat{Y}_{\iota(s)})|^{2}|\mathcal{G}]\di s\right]\notag\\
&&\hspace{-3cm}\leq\ 
k^{2}\E\left[\bigg(|T-t_{n}|^{2}+2|T-t_{n}|\bigg)
\left(\sup_{t\in[0,T\vee t_{n}]}\E\left[|a(\widehat{Y}_{\iota(s)})|^{2} |\mathcal{G}\right]\right)
\right]\notag\\
&&\hspace{-3cm}\leq\ 
k^{2}\left(
\underbrace{
\E\left[\left(|T-t_{n}|^{2}+2|T-t_{n}|\right)^{2}\right]
}_{(\dag)}
\underbrace{
\E\left[
\left(\sup_{t\in[0,T\vee t_{n}]}\E\left[|a(\widehat{Y}_{\iota(t)})|^{2} |\mathcal{G}\right]\right)^{2}
\right]
}_{(\dag\dag)}
\right)^{\frac{1}{2}}\ .\label{aux_CS}
\end{eqnarray}
Note that we have used that $\{t_{i}\}_{i\geq0}$ are measurable with respect to $\mathcal{G}$.
Thanks to part (ii) of Lemmas \ref{overline_integrableT} and \ref{overline_integrableG} we can bound $(\dag\dag)$ by some constant $\kappa_{7}$ depending on $k$ and $T$ only:
\begin{equation*}\label{star}
(\dag\dag)
\leq
\E\left[
\left(\sup_{t\in[0,T]}\E\left[|a(\widehat{Y}_{t})|^{2} |\mathcal{G}\right]
+
\max_{0\leq i \leq n}\E\left[|a(\widehat{Y}_{t_{i}})|^{2} |\mathcal{G}\right]
\right)^{2}
\right]\\
\leq \kappa_{7}\ .
\end{equation*}
To compute the expression in $(\dag)$ we recall that $t_{n}\stackrel{d}{=}\mathbf{g}(n,n/T)$ and hence we apply Jensen's inequality to bound the first three moments of the difference $|T-t_{n}|$ from above by powers of the fourth moment $\E[|T-\mathbf{g}(n,n/T)|^{4}]=3T^{4}(2+n)n^{-3}$, \ie
\begin{eqnarray*}
(\dag)
&=& 
\E[|T-\mathbf{g}(n,n/T)|^{4}]
+
4
\E[|T-\mathbf{g}(n,n/T)|^{3}]
+
4
\E[|T-\mathbf{g}(n,n/T)|^{2}]
\notag\\
&\leq& 
\left(\frac{3T^{4}(2+n)}{n^{3}}\right)
+
4
\left(\frac{3T^{4}(2+n)}{n^{3}}\right)^{3/4}
+
4
\left(\frac{3T^{4}(2+n)}{n^{3}}\right)^{1/2}
\label{star2}\ .
\end{eqnarray*}
Recall now (\ref{aux_CS}) and the upper bounds for $(\dag)$ and $(\dag\dag)$ to conclude the proof. \qed

\subsection{Proof of Theorem \ref{TH_main_g}}

Using the decomposition of the mean square error in (\ref{EP_error}), the proof of the main result of the paper, Theorem \ref{TH_main_g}, which claims that the convergence of the Euler-Poisson scheme towards the solution of the SDE in (\ref{SDE_main}) is of order $\mathcal{O}(n^{-1/2})$ is now a mere corollary of Theorem \ref{unconditional_main} and Proposition \ref{propo_hitting}.

\section{Remarks on the Euler-Poisson scheme}\label{remarksEP}

\subsection{Enhanced Euler-Poisson scheme}\label{EEP_descrip}

The Euler-Poisson scheme has a deterministic number of iterations, but since it is supported on a random grid, the time where the algorithm ends is random. It is straight forward to show that $t_{n}\to T$ a.s. as $n\uparrow\infty$, but it is natural to investigate if there is a more efficient way to stop the algorithm than doing $n$ iterations.

Recall the Poisson process $\mathbf{N}(n/T)$ defined Section \ref{EP_scheme_s} and define $\mathbf{T}(n,T):=t_{\mathbf{N}_{T}+1}$, where we drop the dependence of $n/T$ to ease the notation. Consider the Euler-Poisson scheme now stopped at the random iteration dictated by $\mathbf{N}_{T}+1$. In other words, this enhanced Euler-Poisson scheme considers $\widehat{Y}_{\mathbf{T}(n,T)}$ as the approximation of $Y_{T}$. This modification would be the one that iterates the Euler-Poisson scheme the optimal amount of times to get the final point in the random grid as close as possible to $T$ by overlapping it. Another equivalent modification would be $\widetilde{\mathbf{T}}(n,T):=t_{\mathbf{N}_{T}}$ which makes the final point in the random grid to be the closest one below $T$. A key observation to be made here is that in order to construct $\widehat{Y}_{\mathbf{T}(n,T)}$ we need to be able to sample from the bivariate $(\Delta X_{\mathbf{e}_{i}},\mathbf{e}_{i})$ rather than just from the resolvent of $X$, \ie rather than just from the univariate $\Delta X_{\mathbf{e}_{i}}$. Note that if the distribution of $(\Delta X_{\mathbf{e}_{i}},\mathbf{e}_{i})$ is available then the distribution of $X_{t}$ is given by
\begin{equation}\label{laplace_trns}
\P(\Delta X_{\mathbf{e}(q)}\in\di x, \mathbf{e}(q)\in\di t)= \P(X_{t}\in \di x)\ qe^{-qt}\di t\ .
\end{equation}
In such circumstances one might as well apply the classical Euler scheme for SDEs (also known as the Euler-Maruyama scheme) as it is unlikely to find a process for which $(\Delta X_{\mathbf{e}_{i}},\mathbf{e}_{i})$ is directly available but the distribution of $X_{t}$ is not, in which case the enhanced Euler-Poisson algorithm would avoid the Laplace transformation in (\ref{laplace_trns}). It is an interesting \emph{health check} to derive the rate of convergence for $\E[|Y_{T}-\widehat{Y}_{\mathbf{T}(n,T)}|]$ and see that we recover the asymptotic result of the classical Euler-Maruyama scheme reviewed in the Introduction and suggested by the above argument.  

\begin{cor}\label{C_main_T}
Under the assumptions of Theorem \ref{general_sol_SDE}, we have
\begin{equation*}
\E[|Y_{T}-\widetilde{Y}_{\mathbf{T}(n,T)}|^{2}]\lesssim n^{-1}\ .
\end{equation*}
\end{cor}
\proof We  first proof an analogous result to Proposition \ref{propo_hitting} for the random iteration $\mathbf{N}_{T}+1$. From the construction of $\widehat{Y}$, recall that $\widetilde{Y}_{\mathbf{T}(n,T)}=\widehat{Y}_{\mathbf{T}(n,T)}$, we write
\begin{eqnarray}
\E[|\widehat{Y}_{T}-\widehat{Y}_{\mathbf{T}(n,T)}|^{2}]
&=&\E[|a(\widehat{Y}_{\iota(T)})|^{2}]\E[|X_{\mathbf{T}(n,T)}-X_{T}|^{2}]\notag\\
&\leq&(K_{0}K_{3}+K_{0})\bigg(k^{2}\E[|\mathbf{T}(n,T)-T|^{2}]+2k^{2}\E[|\mathbf{T}(n,T)-T|]\bigg)\notag\\
&=&(K_{0}K_{3}+K_{0})\bigg(k^{2}\frac{T^{2}}{n^{2}}+2k^{2}\frac{T}{n}\bigg)\label{opt_enh}\ ,
\end{eqnarray}
where the only difference with the proof of Proposition \ref{propo_hitting} is the fact 
%we have used (\ref{ineq_a}) and Lemma \ref{overline_integrableT} and \ref{overline_integrableG} to upper bound $a(\widehat{Y}_{\iota(T)})$ and the fact 
that due to the lack of memory property 
$\mathbf{T}(n,T)-T\stackrel{d}{=}\mathbf{e}(n/T)$. To prove the claim of the result we just need to split
the error $|Y_{T}-\widetilde{Y}_{\mathbf{T}(n,T)}|$ into a discretization and a hitting error as done in (\ref{EP_error}) and use the Theorem \ref{unconditional_main} together with (\ref{opt_enh}).\qed

\subsection{Heuristics behind the Euler-Poisson scheme}

The Feynman-Kac representation identifies solutions of certain Partial Integro Differential Equations (PIDE) as conditional expectations of a SDE. This section aims to formalize the relationship between the discretization procedure given by the Euler-Poisson scheme in (\ref{PEs}) and its counterpart in the PIDE representation. We claim that, in some sense, the solution $Y$ of (\ref{SDE_main}) sampled over a random grid generated by the arrival times of a Poisson process is equivalent to perform a discretization in time by the method of lines to the associated Feynman-Kac equation and hence the Euler-Poisson scheme rises as a natural discretization scheme. The relationship is not new and was the basis of Carr \cite{Carr98}, where an approximation for American options of finite maturity is obtained by randomizing the time horizon by an Erlang distribution; Matache et al. \cite{MNS05} also point out informally the relation between a deterministic discretization in time of a  Feynman-Kac PIDE and its probabilistic counterpart. 
\begin{theorem}[{Situ \cite[Section 8.17]{Situ}}]\label{FK_situ}
Consider the following integro-differential operator
\begin{eqnarray*}
\mathcal{A}_{Y}g(x)&:=&
\langle a(x)b,\nabla\rangle g(x)
+
\frac{1}{2}\langle a(x)\Sigma\Sigma^{\text{T}}a^{\text{T}}(x)\nabla,\nabla\rangle g(x)\\
&&\qquad+\ 
\int_{\mathbb{R}^{d_{X}}}\bigg(g(x+a(x)z)-g(x)-\langle a(x)z,\nabla\rangle g(x)\bigg)\Pi(\di z)\ ,
\end{eqnarray*}
taking values in the space $\mathcal{C}^{1,2}([0,T]\times \mathbb{R}^{d_{Y}},\mathbb{R})$
and under the following assumptions
\begin{enumerate}
\item[(i)] $a:=\mathbb{R}^{d_{Y}}\to\mathbb{R}^{d_{Y}}\otimes\mathbb{R}^{d_{X}}$ is bounded \item[(ii)] there exists $\delta_{1}, \delta_{2}>0$ such that $\delta_{1}|\lambda|^{2}\leq\langle a(x)\Sigma\Sigma^{\text{T}}a^{\text{T}}(x)\lambda,\lambda\rangle\leq \delta_{2}|\lambda|^{2}$ for all $x,\lambda\in \mathbb{R}^{d_{Y}}$,
%\item[(iii)] $f:=\mathbb{R}^{d_{Y}}\to\mathbb{R}$ be a given function such that $f\in\mathcal{C}_{0}^{2}$.
%\item[(iv)] $\partial D$ is sufficiently smooth.
\end{enumerate}
together with the assumptions of Theorem \ref{general_sol_SDE}; let $f:=\mathbb{R}^{d_{Y}}\to\mathbb{R}$ be a bounded continuous function, \ie $f\in\mathcal{C}_{0}$. If $u(t,x)\in \mathcal{C}^{1,2}([0,T]\times \mathbb{R}^{d_{Y}},\mathbb{R})$ is a classical solution of the PIDE
\begin{equation}\label{FK_PDE}
\frac{\partial}{\partial t}u(t,x)=\mathcal{A}_{Y}u(t,x)
\end{equation}
with initial condition $u(0,x)=f(x)$, then
\begin{equation}\label{cond_exp}
u(T-t,x)=\E[f(Y_{T})|Y_{t}=x]=\E[f(Y_{T-t})|Y_{0}=x]:=\E_{x}[f(Y_{T-t})]\ ,
\end{equation}
where $Y$ is the unique strong solution of (\ref{SDE_main}) and $0\leq t\leq T$.
\end{theorem}
The converse of the preceding statement also holds with appropriate assumptions. The above result can be written under much more general assumptions and in terms of weak solutions of the PIDE, but the simpler statement above is enough to make the point in this section. A typical setting where the above relation is exploited happens when (\ref{cond_exp}) represents the price of an option under the risky asset $Y$ which is computed by numerically solving the associated PIDE. The celebrated Black-Scholes formula is an example of this approach when the underlying follows a geometric Brownian motion; for incomplete markets generated by L\'evy processes similar formulas hold (cf. Chan \cite{Chan99}). 

Recall the random times $\{t_{i}\}_{i\geq0}$ defined in Section \ref{EP_scheme_s} as the arrival times of a Poisson process $\mathbf{N}$, and consider the Laplace-Carlson transform, $\mathcal{L}$, of $u(t,x)$, that is
\begin{eqnarray}
\mathcal{L}[u](x)&:=&\int_{0}^{\infty}\frac{n}{T}\exp(-nt/T)u(t,x)\di t\notag\\
&=&
\int_{0}^{\infty}\frac{n}{T}\exp(-nt/T)\E_{x}[f(Y_{t})]\di t\notag\\
&=&\E_{x}[f(Y_{\mathbf{e}(n/T)})]=\E_{x}[f(Y_{t_{1}})]\label{random_1}\ ,
\end{eqnarray}
where we have used the boundedness of $f\in\mathcal{C}_{0}$ to apply Fubini's theorem. Note that the last term in the above equation corresponds to the expectation of the solution in (\ref{SDE_main}) at the first arrival time of the Poisson process $\mathbf{N}$. Moreover, due to the boundedness of $f$ we can also interchange the differential operator $\mathcal{A}_{Y}$ and the transform $\mathcal{L}$ to obtain the non-homogeneous integro-differential equation satisfied by the Laplace-Carlson transform: 
\begin{equation}\label{LC_PDE}
\frac{\mathcal{L}[u](x)-f(x)}{T/n}=\mathcal{A}_{Y}\mathcal{L}[u](x)\ .
\end{equation}
Indeed, the above equation transforms the differential $\frac{\partial}{\partial t}$ in (\ref{FK_PDE}) into a difference which turns out to be of the same form as the first order finite difference approximation in time of (\ref{FK_PDE}) with respect to $\mathcal{L}[u]$ instead of $u$ due to the homogeneity of $\mathcal{A}_{Y}$. To fix ideas, the following proposition explicitly relates the solution $Y$ at the arrival times of $\mathbf{N}$ with the iterates of what is known as the method of lines or Rothe's method in the literature of numerical analysis of PIDEs.
\begin{prop}\label{PDE_exponential}
Under the assumptions of Theorem \ref{general_sol_SDE} and \ref{FK_situ}, consider the Rothe's discretization of (\ref{FK_PDE}) given by 
\begin{equation}\label{Roths}
\frac{u_{i}(x)-u_{i-1}(x)}{T/n}=\mathcal{A}_{Y}u_{i}(x)\ ,
\end{equation}
for $i=1,\ldots,n$ with $u_{0}(x)=f(x)$. Then, for all $i=1,\ldots,n$, 
\begin{equation*}
u_{i}(x)=\E_{x}[f(Y_{t_{i}})]\ .
\end{equation*}
\end{prop}
\proof
It is clear that the solution of (\ref{SDE_main}) given by Theorem \ref{general_sol_SDE} has the strong Markov property (cf. Protter \cite[Theorem 32 p. 294]{P05}). Therefore, we write
\begin{equation*}
\E_{x}[f({Y}_{t_{i}})]
=\E_{x}[\E_{{Y}_{t_{1}}}[\E_{{Y}_{t_{2}}}[\cdots\E_{{Y}_{t_{i-1}}}[f(Y_{t_{i}})]\cdots]]]\ ,
\end{equation*}
and apply recursively the arguments derived from (\ref{random_1}) and (\ref{LC_PDE}) in the above nested expectations to obtain the recursive solutions that solve the system of differential equations in (\ref{Roths}).\qed

The justification of the Euler-Poisson scheme was already discussed in Section \ref{main_result_discussion} through the description of the scenarios where the method is a feasible option. The above result strengthens that discussion as Proposition \ref{PDE_exponential} is identifying the Euler-Poisson scheme to an approximation of the classical Rothe's method.

\subsection{Pathwise convergence}

The Euler-Poisson scheme is supported on a random grid and there is no straightforward way to perform a pathwise numerical analysis of the algorithm. Nevertheless the above analogy with Rothe's discretization method suggest that one may try to study the behavior of
\begin{equation*}
\E[
\max_{1\leq i\leq n}
|Y_{iT/n}-\widetilde{Y}_{t_{i}}|^{2}
]=\E[
\max_{1\leq i\leq n}
|Y_{iT/n}-\widehat{Y}_{t_{i}}|^{2}
]\ .
\end{equation*}
Indeed, Theorem \ref{unconditional_main} states a pathwise result for the discretization error and hence, using the decomposition in (\ref{EP_error}), one would only need to obtain a pathwise analogue of the hitting error in order to study the above quantity, \ie a pathwise generalization of Proposition \ref{propo_hitting}. Unfortunately the latter is not true. 

A weaker statement which can be proved, involving the entire path of the Euler-Poisson scheme, is
\begin{equation}\label{max_out}
\max_{1\leq i\leq n}\E[|Y_{iT/n}-\widetilde{Y}_{t_{i}}|^{2}]\leq K_{8} n^{-1/2}\ ,
\end{equation}
where $K_{8}$ is a positive constant depending on $k$ and $T$ only. The result in (\ref{max_out}) is a direct consequence of Theorem \ref{unconditional_main} and Proposition \ref{propo_hitting} together with the following lemma:
\begin{lemma}\label{max_grid}
The following bound holds
\begin{equation*}
\E\left[\max_{1\leq i\leq n}\left|t_{i}-\frac{Ti}{n}\right|^{p}\right]\leq8
\E\left[\left|t_{n}-T\right|^{p}\right]\qquad\text{for}\quad p\geq1\ .
\end{equation*}
\end{lemma}
\proof Define $Z_{i}:=\mathbf{e}_{i}(n/T)-\frac{T}{n}$ for $1\leq i\leq n$ and $Z_{i}:=0$ for $i>n$. The sequence $\{Z_{i}\}_{i\geq1}$ is centred and mutually independent and the result follows from Doob \cite[Theorem 5.1]{Doob53}.\qed

\appendix
\section{Moments of $\tau$}\label{Appendix}

It is well understood that given the Poisson process $\mathbf{N}$ conditioned on having $m$ arrivals up to time $T$, the location along $[0,T]$ of those $m$ arrival times have the same distribution as $m$ ordered independent uniform random variables on $[0,T]$. Therefore in order to study the random variable $\tau$ defined in (\ref{def_tau}) we can start by studying the largest partition on the interval $[0,1]$ defined by $m$ independent uniform random variables in $[0,1]$. 

Given $m>0$, let $\{U_{i}\}_{i=1,\ldots,m-1}$ be a sequence of i.i.d. random variables with common uniform distribution in $[0,1]$ and consider $U_{(i)}$ for $i=0,\ldots,m$ its order statistics where $U_{0}=0$ and $U_{m}=1$. Denote the largest gap by 
\begin{equation*}
\lambda_{m}:=\max_{i=1\ldots,m}\{U_{(i)}-U_{(i-1)}\}\ .
\end{equation*}
Recall the definition of $\tau$ in (\ref{def_tau}).
Then the conditional distribution of $\tau$ is, up to a constant, equal to $\lambda$. Indeed, 
$\frac{1}{T}\tau$ is, conditioned on $\mathbf{N}_{T}$, equal in distribution to $\lambda_{\mathbf{N}_{T}+1}$. In particular we have 
\begin{equation}\label{cond_tau}
\frac{1}{T}\E[\tau]=\E[\lambda_{\mathbf{N}_{T}+1}]\ .
\end{equation}
Fisher \cite{Fisher29} already studied the behaviour of $\lambda_{m}$ and the following expression can be found in Mauldon \cite{Mauldon51}:
\begin{equation*}
\E[(1-\lambda_{m} s)^{-m}]=\frac{m!}{1-s}\prod_{j=2}^{m}\frac{1}{j-s}\qquad |s|<1/2\ ,m\geq1\ .
\end{equation*}
All moments of $\lambda_{m}$ can be expanded form the above expression and in particular, for $m\geq1$, we have
\begin{equation*}
\E[\lambda_{m}]=\frac{\sum_{j=1}^{m}\frac{1}{j}}{m}=\frac{\Psi(m+1)+\gamma}{m}\ ,
\end{equation*}
where $\Psi$ is the digamma function (see Abramowitz and Stegun \cite[Section 6.3.2 and Section 6.4.10]{AS70}). Recall that the function $\Psi(m+1)+\gamma$ is zero for $m=0$, positive for $m>0$ and grows asymptotically as $\log(m+1)$, \ie $\lim_{m\to\infty}\Psi(m)/\log(m)=1$. Therefore, there is a constant $\kappa_{0} >0$ independent of $m$ such that $\Psi(m+1)+\gamma\leq \kappa_{0}\log(m+1)$ for  $m\geq1$,
hence 
\begin{equation*}
\E[\lambda_{m}]\leq \kappa_{0}\frac{\log(m+1)}{m}
\qquad
\text{ for }
\qquad m=1,2,\ldots
\end{equation*}
\begin{prop}\label{tau_Poisson}
Under the above notations, 
\begin{equation*}
\E[\tau]+\E[\tau^{2}]\lesssim n^{-1}\log(n)\ .
\end{equation*}
\end{prop}
\proof
According to (\ref{cond_tau}) and recalling that the arrival rate for $\mathbf{N}_{T}$ is $n/T$, we have
\begin{eqnarray*}
\frac{1}{T}\E[\tau]=\E[\lambda_{\mathbf{N}_{T}+1}]
&=&\sum_{k=0}^{\infty}\E[\lambda_{\mathbf{N}_{T}+1}|\mathbf{N}_{T}=k]\P(\mathbf{N}_{T}=k)\\
&\leq&\kappa_{0}\sum_{k=0}^{\infty}\frac{\log(k+2)}{k+1}\exp\left(-\frac{n}{T}\right)\frac{\left(\frac{n}{T}\right)^{k}}{k!}\\
&=&\kappa_{0}\frac{T}{n}\sum_{k=1}^{\infty}\log(k+1)\exp\left(-\frac{n}{T}\right)\frac{\left(\frac{n}{T}\right)^{k}}{k!}\\
&=&\frac{\kappa_{0} T}{n}\E[\log(\mathbf{N}_{T}+1)]\leq\frac{\kappa_{0} T}{n}\log(\E[\mathbf{N}_{T}]+1)=\frac{\kappa_{0} T\log\left(\frac{n}{T}+1\right)}{n}
\end{eqnarray*}
where the last inequality is due to Jensen's inequality and the concavity of $x\to\ln(x+1)$ for $x\geq0$. To derive the claim of the proposition we can crudely upper bound $\lambda_{m}^{2}\leq\lambda_{m}$, since $\lambda_{m}\in[0,1]$, and hence
\begin{equation*}
\frac{1}{T^{2}}\E[\tau^{2}]=\E[\lambda_{\mathbf{N}_{T}+1}^{2}]\leq\E[\lambda_{\mathbf{N}_{T}+1}]\ .
\end{equation*}
\qed

\end{document}